\documentclass[conference,a4paper]{IEEEtran}

\usepackage{bm,bbm}
\usepackage{amsmath,amssymb}
\usepackage{graphicx,subfigure}
\usepackage{units}
\usepackage{url}
\usepackage{cite}
\usepackage{booktabs}

\begin{document}

\title{Contrast Measures based on the Complex Correlation Coefficient for PolSAR Imagery}
\author{
\IEEEauthorblockN{Alejandro C.\ Frery$^1$, Renato J.\  Cintra$^2$, Abra\~ao D.\ C.\ Nascimento$^3$}
\IEEEauthorblockA{
$^1$LaCCAN -- Laborat\'orio de Computa\c c\~ao Cient\'ifica e An\'alise Num\'erica\\
Universidade Federal de Alagoas\\
Av. Lourival Melo Mota, s/n\\
57072-900 Macei\'o -- AL, Brazil\\
\texttt{acfrery@gmail.com}\\
$^2$Departamento de Estat\'istica\\
Universidade Federal de Pernambuco\\
50740-540 Recife -- PE, Brazil\\
\texttt{rjdsc@de.ufpe.br}\\
$^3$Departamento de Estat\'istica\\
Universidade Federal de Para\'iba\\
58051-900 Jo\~ao Pessoa -- PB, Brazil\\
\texttt{abraao@de.ufpb.br}
} 
}

\maketitle

\begin{abstract}
We derive contrast measures which involve the number of looks and the complex correlation coefficient between polarization channels in PolSAR imagery.
Using asymptotic results which characterize the behavior of these measures, we derive statistical regions of confidence which lead to test of hypothesis.
An application to real data is performed, confirming the importance of the proposals.
\end{abstract}

\section{Introduction}

The aim of remote sensing is to capture and to analyze information scenes concerning the Earth surface.
In this context, polarimetric synthetic aperture radar (PolSAR) has achieved a prominent position among the remote sensing technologies~\cite{RemoteSensingPolarimetricRadarMott}.
Such systems employ \textit{coherent illumination} in its processing and, as a consequence, their resulting images are contaminated with fluctuations on its detected intensity called ``speckle''.
Speckle significantly degrades the perceived image quality, as much as the ability of extracting information from the data.

Speckle is well described by statistical models.
Thus, two pre-processing steps are often sought:
(i)~the identification of a probability distribution for PolSAR image regions~\cite{FreitasFreryCorreia:Environmetrics:03}, and 
(ii)~the derivation of statistical methods for quantifying contrast between such regions~\cite{FreryCintraNascimento2013}.

A successful statistical model for homogeneous regions in PolSAR images is the scaled complex Wishart law~\cite{FreitasFreryCorreia:Environmetrics:03}.
This distribution is equipped with two parameters: the number of looks and the complex covariance matrix.
A rich discussion about estimation and interpretation of the number of looks was given by Anfinsen~\textit{et al.}~\cite{EstimationEquivalentNumberLooksSAR}.
In terms of the covariance matrix, it contains all necessary information  to characterize the backscattered data.
 Conradsen~\textit{et al.}~\cite{Conradsen2003} discuss hypothesis tests based on the covariance matrix.

Lee~\textit{et al.}~\cite{leeetal1994b} proposed a reparametrization for the covariance matrix in terms of the complex correlation coefficient.
Statistical models which describe the multilook phase difference, the magnitude of the complex product, and both intensity and amplitude ratios between two components of the scattering matrix are provided in that work.
The resulting density functions have closed forms which depend on the complex correlation coefficient and on the number of looks.

Many authors have utilized the complex correlation coefficient as an important quantity for analyzing PolSAR images. 
For instance, Ainsworth~\textit{et al.}~\cite{Ainsworth20082876} presented evidence that the complex correlation coefficient between channels can be utilized to identify man-made targets.

In recent years, information-theoretic based measures have used to derive new PolSAR image processing methods.
Literature in this field of research can be divided in two groups:
(i)~works that involve deterministic tools~\cite{MorioRefregierGoudailFernandezDupuis2009},
and 
(ii)~contributions that consider statistical properties of these measures~\cite{FreryCintraNascimento2012,FreryCintraNascimento2013}.
In this work, we  advance the statistical inference based on the complex correlation coefficient in for PolSAR data --- a contribution within the second group of works. 

In summary, our contributions are two-fold:
\begin{enumerate}
\item Based on the parametrization proposed by Lee~\textit{et al.}~\cite{leeetal1994b}, we derive four contrast measures which depend on the complex correlation coefficient.	
These measures were obtained considering four distances from the $h$-$\phi$ class of distances proposed by Salicr\'u et al.~\cite{salicruetal1994}: the Kullback-Leibler, R\'enyi (of order $\beta$), Bhattacharyya, and Hellinger distances between reparametrized scaled complex Wishart distributions.

\item We study the asymptotic properties of these measures, and we propose new confidence regions for these distances which allow comparing two PolSAR regions.
\end{enumerate}

\section{The model}

If the complex return with $p$ polarization channels follows a complex Gaussian law~\cite{FreitasFreryCorreia:Environmetrics:03}, the multilook covariance matrix return  $\boldsymbol{Z}$ follows a scaled complex Wishart distribution characterized by the following probability density function:
\begin{equation}
 f_{\boldsymbol{Z}}({Z};\boldsymbol{\Sigma},L) = \frac{L^{pL}|{Z}|^{L-p}}{|\boldsymbol{\Sigma}|^L \Gamma_p(L)} \exp\bigl(
-L\operatorname{tr}\bigl(\boldsymbol{\Sigma}^{-1}{Z}\bigr)\bigr),
\label{eq:denswishart}
\end{equation}
where $\Gamma_p(L)=\pi^{p(p-1)/2}\prod_{i=0}^{p-1}\Gamma(L-i)$, $L\geq p$, $\Gamma(\cdot)$ is the gamma function, 
$\operatorname{tr}(\cdot)$ is the trace operator,
and $\exp(\cdot)$ is the exponential function.
The parameters which index this distribution are the number of looks $L$, and the covariance matrix $\boldsymbol{\Sigma}$.
This situation is denoted $\boldsymbol{Z}\sim \mathcal W(\boldsymbol{\Sigma},L)$, and this distribution satisfies $\operatorname{E}\{\boldsymbol{Z}\}=\boldsymbol{\Sigma}$, which is a Hermitian positive definite matrix~\cite{EstimationEquivalentNumberLooksSAR}.

Lee \textit{et al}.~\cite{leeetal1994b} presented a reparametrization of a particular case of the complex Wishart distribution based on a two-element scattering vector $\boldsymbol{y}^{(k)}=\big[y_1^{(k)},\,y_2^{(k)}\big]^\top$ at the $k$th look, for which the covariance matrix is written as
$$
 \boldsymbol{Z}= \left[\begin{array}{cc}
 {z}_{11}& \alpha e^{\mathbf{i}\Delta}\\
\alpha e^{-\mathbf{i}\Delta} &  {z}_{22}
\end{array} \right],
$$
where $(\cdot)^\top$ is the transposition operator, 
$\alpha$ and $\Delta$ are the sample multilook magnitude and phase, respectively, 
and $z_{ii}=L^{-1}\sum_{k=1}^L y_i^{(k)} y_i^{(k)*}$, for $i=1,2$.
The resulting covariance matrix is 
\begin{equation}
\boldsymbol{\Sigma}=\left[\begin{array}{cc}
{\sigma}_{11}& \sqrt{{\sigma}_{11}{\sigma}_{22}} |{\rho}_c| e^{\mathbf{i}\delta}\\
\sqrt{{\sigma}_{11}{\sigma}_{22}} |{\rho}_c| e^{-\mathbf{i}\delta} & {\sigma}_{22}
\end{array} \right],
\label{matcov}
\end{equation}
where ${\sigma}_{ii}=\operatorname{E}(z_{ii})$, $\delta$ is the population multilook phase, and ${\rho}_c$ is the complex correlation coefficient between ${z}_{11}$ and ${z}_{22}$.

In order to study the preservation of polarimetric properties along the process of filtering PolSAR imagery, Lee \textit{et al.}~\cite{LeeandGrunes1999} studied the correlation coefficient between polarization channels.
It is given by
\begin{equation}
\rho_c=\frac{\operatorname{E}(S_{\text{HH}}S_{\text{VV}}^*)}{\sqrt{\operatorname{E}(|S_{\text{HH}}|^2)\operatorname{E}(|S_{\text{VV}}|^2)}}.
\label{cor}
\end{equation}
The normalized quantities
$$
B_1=\frac{{z}_{11}}{{\sigma}_{11}}, \quad B_2=\frac{ {z}_{22}}{{\sigma}_{22}},\quad \eta=\frac{\alpha}{\sqrt{{\sigma}_{11}{\sigma}_{22}}},
$$
and $\Delta$ obey the distribution characterized by the following joint probability density function:
\begin{align}
f(B_1,B_2,\eta,\Delta;\rho_c,L)=\frac{\eta(B_1B_2-\eta^2)^{L-2}}{\pi(1-|{\rho}_c|^2)^L \Gamma{(L)}\Gamma{(L-1)}} \nonumber \\
\times \exp \Big\{-\frac{B_1+B_2-2\eta |{\rho}_c| \cos(\Delta-\delta)}{1-|{\rho}_c|^2} \Big\}.
\label{correlationinjoint}
\end{align}

The complex correlation coefficient has been utilized as a important quantity for identifying contrast in PolSAR images.
For instance, Lee \textit{et al.}~\cite{LeeandGrunes1999} presented results which provide evidences that 
changes of correlation coefficient between polarization channels can be captured when one considers pixels of different regions.

\section{Stochastic Distances based on the Complex Correlation Coefficient}\label{Sec:Distances}

We adhere to the convention that a (stochastic) ``divergence'' is any non-negative function between two probability measures which obeys the identity of definiteness property~\cite{FreryCintraNascimento2013}.
If the function is also symmetric, it is called ``distance''~\cite[ch. 1 and 14]{deza2009encyclopedia}.

An image can be understood as a set of regions, in which the enclosed pixels are observations of random variables following a certain distribution.
Therefore, stochastic dissimilarity measures can be used to assess the difference between the distributions that describe different image areas~\cite{HypothesisTestingSpeckledDataStochasticDistances}.

Dissimilarity measures were submitted to a systematic and comprehensive treatment in~\cite{salicruetal1994}, leading to the proposal of the class of $(h,\phi)$-divergences.
Stochastic distances applied to intensity SAR data were presented in~\cite{HypothesisTestingSpeckledDataStochasticDistances,ParametricNonparametricTestsSpeckledImagery} and to PolSAR models in~\cite{FreryCintraNascimento2013}.

We consider here four stochastic distances between the models characterized by the densities $f_1$ and $f_2$ with the same support $\mathcal A$:
\begin{itemize}
\item [i)] Kullback-Leibler:
$
d_{\text{KL}} = \frac 12\int_{\mathcal A}(f_1-f_2)
\log\frac{f_1}{f_2}.
$
\item [ii)] R\'{e}nyi of order $\beta\in(0,1)$:
$ 
d_ {\text{R}}^{\beta} = (\beta-1)^{-1}
 \log( \int_{\mathcal A} f_1^{\beta} f_2^{1-\beta}
+ \int_{\mathcal A} f_1^{1-\beta} f_2^{\beta})/2.
$
\item [iii)] Bhattacharyya:
$
d_{\text{BA}} = -\log\int_{\mathcal A}\sqrt{f_1 f_2}.
$
\item [iv)] Hellinger:
$
d_{\text{H}} = 1-\int_{\mathcal A}\sqrt{f_1 f_2}.
$
\end{itemize}

In practical applications, the densities $f_1$ and $f_2$ are not know.
Their parameters are usually estimated with samples of sizes $N_1$ and $N_2$, yielding $\widehat{f_1}=f_1(\cdot;\widehat{\boldsymbol{\theta_1}}(N_1))$ and $\widehat{f_2}=f_2(\cdot;\widehat{\boldsymbol{\theta_2}}(N_2))$.
Whenever there is no risk of ambiguity, i.e., when the distribution is the same, only the estimators can be used to denote the distances.

These distances become more useful and comparable scaling them into test statistics, as discussed in~\cite{HypothesisTestingSpeckledDataStochasticDistances}:
\begin{equation}
S_{\mathcal D}\big(\widehat{\boldsymbol{\theta_1}}{(N_1)},\widehat{\boldsymbol{\theta_2}}{(N_2)}\big) = \frac{2 N_1 N_2 v_{\mathcal D}}{N_1+N_2}
d_{\mathcal D}\big(\widehat{\boldsymbol{\theta_1}}{(N_1)},\widehat{\boldsymbol{\theta_2}}{(N_2)}\big),
\label{eq:TestStatistic}
\end{equation}
where $v_{\mathcal D} = 1,\,\beta^{-1},\, 4,\text{ and } 4$ for ${\mathcal D}=$ KL, R, BA, and H, respectively, and $\widehat{\boldsymbol{\theta_1}}{(N_1)}=[\widehat{\boldsymbol{\Sigma}}(N_1),\widehat{L}(N_1)]$ and $\widehat{\boldsymbol{\theta_2}}{(N_2)}=[\widehat{\boldsymbol{\Sigma}}(N_2),\widehat{L}(N_2)]$ are the maximum likelihood estimators for $\boldsymbol{\theta_1}$ and $\boldsymbol{\theta_2}$ using different random samples of sizes $N_1$ and $N_2$, respectively.
Under mild conditions, $S_{\mathcal D}\big(\widehat{\boldsymbol{\theta_1}}{(N_1)},\widehat{\boldsymbol{\theta_2}}{(N_2)}\big)$ is asymptotically distributed as a $\chi^2_M$ random variable, where $M$ is the dimension of the parameter $\boldsymbol{\theta}$.

\section{Results}

This section presents contrast measures based on the number of looks and the complex correlation coefficient,
and provides probabilistic criteria for discriminating two PolSAR regions in terms of their sample correlation coefficients between polarization channels.

Since this seems untractable for random vectors equipped with densities given by Eq.~\eqref{correlationinjoint}, we derive the four distances discussed in the previous section between the random matrices $\bm{Z}_1$ and $\bm{Z}_2$ such that
$$
\bm{Z}_k=\left[\begin{array}{cc}
 {z}_{11}^{(k)}& \alpha e^{\mathbf{i}\Delta^{(k)}}\\
\alpha e^{-\mathbf{i}\Delta^{(k)}} &  {z}_{22}^{(k)}
\end{array} \right]
$$
and
$$
\operatorname{E}(\bm{Z}_k)=\bm{\Sigma}_k=
\left[\begin{array}{cc}
{\sigma}_{11}^{(k)}& \sqrt{{\sigma}_{11}^{(k)}{\sigma}_{22}^{(k)}} |{\rho}_c| e^{\mathbf{i}\delta}\\
\sqrt{{\sigma}_{11}^{(k)}{\sigma}_{22}^{(k)}} |{\rho}_c| e^{-\mathbf{i}\delta} & {\sigma}_{22}^{(k)}
\end{array} \right],
$$
for $k=1,2$.
In order to obtain
closed-form contrast measures, 
we assume that $\Pr(\Delta^{(k)}=0)=1$ 
or, equivalently, 
that the population multilook phase is zero.
Thus, the used distances are denoted by
\begin{equation}
d_{\mathcal D}(\boldsymbol{\theta}_1,\boldsymbol{\theta}_2)
\equiv
d_{\mathcal D}(\bm{Z}_1,\bm{Z}_2),
\label{generalidentity}
\end{equation}
where $\boldsymbol{\theta}_k=(\rho_k,L_k,\sigma_1^{(k)},\sigma_2^{(k)})$, 
for $k=1,2$.

A common practice when dealing with observations that have different scales or units is to standardize the variables~\cite{mardia}. 
Such transformation may be relevant since the mean intensity in the cross-polarized channels can be quite different from the one in the co-polarized channels~\cite{Ainsworth20082876}. 
When that is the case, only the correlation coefficients are left to check if two samples come from the same distribution.
We, thus, derived the distances $d_{\mathcal D}(\rho_1,\rho_2\mid L)\equiv d_{\mathcal D}([\rho_1,L]^\top,[\rho_2,L]^\top)$, when $\sigma_{11}=\sigma_{22}=1$ and  $L_1=L_2=L$ in~\eqref{generalidentity}.
They are given by Eqs.~\eqref{dklpho}-\eqref{dHpho}.

Note that the change of correlation coefficient between polarization channels 
can be captured 
when one considers pixels of different regions.
In the subsequent discussion, 
we provide two methodologies in terms of Eqs.~\eqref{dklpho} and~\eqref{dHpho} which 
aim to discriminate two PolSAR regions
based on the difference between $\rho_1$ and $\rho_2$.
%

From the asymptotic result discussed in the last paragraph of Section~\ref{Sec:Distances},
four new hypothesis tests 
are derived for checking whether two samples, of sizes $N_1$ and $N_2$, 
come from regions with statistically similar correlation coefficients between channels.
As a consequence, 
the resulting statistics can be also used as confidence regions.
For simplicity, 
we present only the results relative to the Kullback-Leibler and Hellinger distances.

The statistics based on Kullback-Leibler and Hellinger distances present the following confidence regions,
respectively:
Given a PolSAR image, 
let $L$ be the number of looks known and constant on this image,
\begin{align}
&R_\text{KL}(\rho_1,\rho_2\mid L)=\Big\{(\rho_1,\rho_2)\in \mathbb{C}\times\mathbb{C}\colon\nonumber \\
&\underbrace{\frac{(1-|\rho_1||\rho_2|)(2-|\rho_1|^2-|\rho_2|^2)}{(1-|\rho_1|^2)(1-|\rho_2|^2)}}_{\xi_1}
\leq \underbrace{\frac{(N_1+N_2)}{2\,L\,N_1\,N_2}\chi^2_1(\eta) +2 }_{t_\text{KL}}\Big\}
\label{cregion1}
\end{align}
and
\begin{align}
&R_\text{H}(\rho_1,\rho_2\mid L)=\Big\{(\rho_1,\rho_2)\in \mathbb{C}\times\mathbb{C}\colon \nonumber \\
&\underbrace{\frac{\sqrt{(1-|\rho_1|^2)(1-|\rho_2|^2)}}{4-(|\rho_1|+|\rho_2|)^2}}_{\xi_2}
\geq \underbrace{\frac 14 \Big[1- \frac{(N_1+N_2)}{2\,L\,N_1\,N_2}\,\chi^2_1(\eta)\Big]^{1/L}}_{t_\text{H}} \Big\},
\label{cregion2}
\end{align}
where $\eta$ is the specified nominal level.

Given two polarimetric data samples, 
their correlation coefficients can be estimated according to Eq.~\eqref{cor}.
Quantities  $\widehat\xi_1$, $\widehat\xi_2$, $t_\text{KL}$, and $t_\text{H}$ 
are evaluated by means of \eqref{cregion1} and \eqref{cregion2},
replacing $\rho_1$ and $\rho_2$ by their  estimates (or sample counterparts) $\widehat \rho_1$ and $\widehat \rho_2$. 
We then propose the following decision rules:
\begin{itemize}
\item Kullback-Leibler criterion: If $\widehat\xi_1  \leq  t_\text{KL}$, then we have statistical evidence that the samples come from populations with similar correlation between HH and VV channels.
\item Hellinger criterion: If $\widehat\xi_2 \geq t_\text{H}$, the samples present equivalent correlation between channels.
\end{itemize}
Thus, regions $R_\text{KL}(\rho_1,\rho_2\mid L)$ and $R_\text{H}(\rho_1,\rho_2\mid L)$ along with~\eqref{cor} are methods for comparing two regions based on their estimates for the correlation coefficient. 

\begin{figure}[htb]
\centering
\includegraphics[width=1\linewidth]{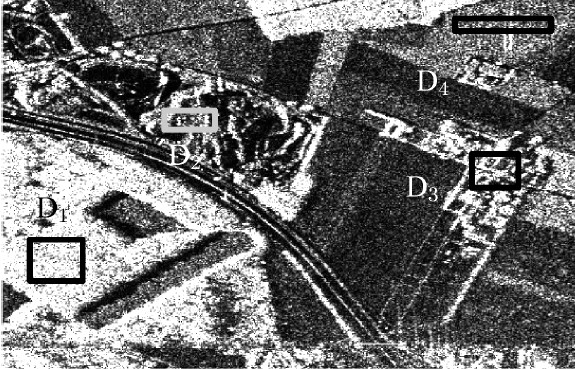}
\caption{PolSAR image (HH channel) with four samples.}
\label{correlationfigu}
\end{figure}

In order to illustrate this methodology, consider the samples highlighted in Fig.~\ref{correlationfigu} (which is extracted from an E-SAR image from surroundings of We\ss ling, Germany).
Table~\ref{correlation} lists the quantities observed, along with the decisions they led to.
Notice that both rules discriminate well.

\begin{table}[htb]                                               
\centering                                                        
\caption{Results based on confidence regions}  
\label{correlation}                                               
\begin{tabular}{c@{\quad}c@{\,\,}c@{\quad}c@{\quad}c@{\,\,}c@{\quad}c@{}}\toprule                                    
 Regions & $\widehat\xi_{1}$ &  $t_\text{KL}$  & Decision & $\widehat\xi_{2}$  & $t_\text{H}$ & Decision\\ \midrule
$\text{D}_1$-$\text{D}_2$ & 2.0220 & 2.0016 & \textit{Distinct} & 0.2486 & 0.2499 & \textit{Distinct}\\       
$\text{D}_1$-$\text{D}_3$ & 2.0068 & 2.0012 & \textit{Distinct} & 0.2495 & 0.2499 & \textit{Distinct}\\       
$\text{D}_1$-$\text{D}_4$ & 2.0734 & 2.0014 & \textit{Distinct} & 0.2455 & 0.2499 & \textit{Distinct}\\       
$\text{D}_2$-$\text{D}_3$ & 2.0534 & 2.0018 & \textit{Distinct} & 0.2467 & 0.2499 & \textit{Distinct}\\       
$\text{D}_2$-$\text{D}_4$ & 2.0149 & 2.0020 & \textit{Distinct} & 0.2490 & 0.2498 & \textit{Distinct}\\       
$\text{D}_3$-$\text{D}_4$ & 2.1254 & 2.0016 & \textit{Distinct} & 0.2424 & 0.2499 & \textit{Distinct}\\       
\bottomrule
\end{tabular}                                                     
\end{table}

\begin{figure*}[tb]
\begin{equation}
d_\text{KL}(\rho_1,\rho_2\mid L)=L\bigg[\frac{(1-|\rho_1||\rho_2|)(2-|\rho_1|^2-|\rho_2|^2)}{(1-|\rho_1|^2)(1-|\rho_2|^2)}-2\bigg]
\label{dklpho}
\end{equation}
\begin{equation}
d_\text{R}^{\beta}(\rho_1,\rho_2\mid L)=\frac{\log 2}{1-\beta}+ \frac{1}{\beta-1}
\log\Bigg\{
\Bigg[\frac{(1-|\rho_1|^2)^{1-\beta}(1-|\rho_2|^2)^{\beta}}{1-\{|\rho_1|-\beta(|\rho_1|-|\rho_2|)\}^2}\Bigg]^L 
+\Bigg[\frac{(1-|\rho_2|^2)^{1-\beta}(1-|\rho_1|^2)^{\beta}}{1-\{|\rho_2|-\beta(|\rho_2|-|\rho_1|)\}^2}\Bigg]^L\Bigg\}
\label{drepho}
\end{equation}
\begin{equation}
d_\text{B}(\rho_1,\rho_2\mid L)=L\bigg\{\frac{\log(1-|\rho_1|^2)+\log(1-|\rho_2|^2)}{2} 
-2\log 2+\log\bigg[\frac{4-(|\rho_1|+|\rho_2|)^2}{(1-|\rho_1|^2)(1-|\rho_2|^2)}\bigg]\bigg\}
\label{dBpho}
\end{equation}
\begin{equation}
d_\text{H}(\rho_1,\rho_2\mid L) =1-\Bigg[4\frac{\sqrt{(1-|\rho_1|^2)(1-|\rho_2|^2)}}{4-(|\rho_1|+|\rho_2|)^2}\Bigg]^L
\label{dHpho}
\end{equation}
\hrulefill
\end{figure*}
  
\begin{figure}[htb]
\centering
\subfigure[$L=2$\label{phonew10}]{\includegraphics[width=.48\linewidth]{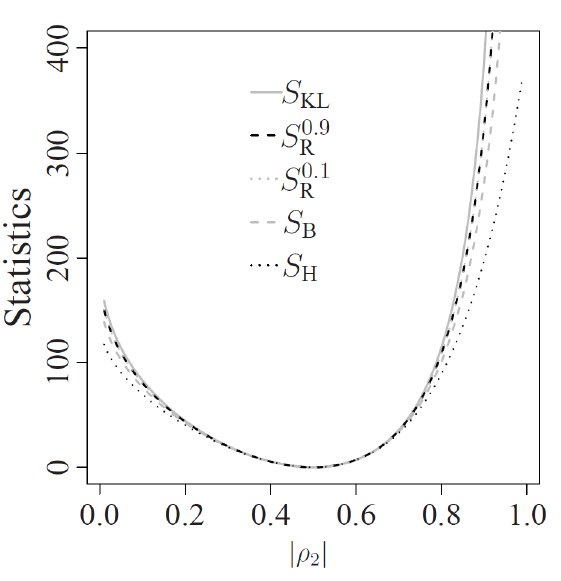}}
\subfigure[$L=5$\label{phonew20}]{\includegraphics[width=.48\linewidth]{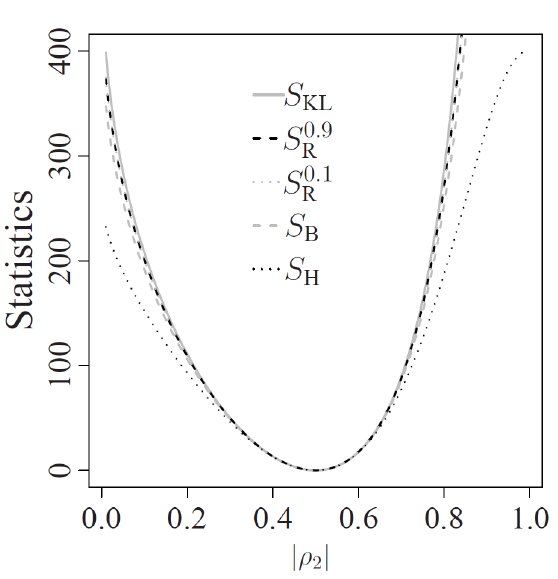}}
\caption{Sensitivity of statistics in~\eqref{dklpho}-\eqref{dHpho} for $L_1=L_2=L=\{2,5\}$, and $|\rho_1|^2=0.5$.}
\label{phocurves}
\end{figure}

Statistics derived from these distances can be used to measure the influence of the number of looks on the complex correlation coefficient in $\boldsymbol{Z}$.   
Fig.~\ref{phocurves} illustrates the resulting statistics for $|\rho_1|^2=0.5$ and $L_1=L_2=L\in\{2,5\}$ (Figures~\ref{phonew10} and~\ref{phonew20}, respectively), while $|\rho_2|$ varies on $[0,1]$.
As expected, the curves have their minimum value, zero, at $\rho_1=\rho_2$.
Notice, however, that the curves are steeper on the interval $(|\rho_1|,1)$ than on $(|\rho_1|,0)$.
Thus for discrimination purposes, the statistics provide better sensitivity capabilities when $|\rho_2|>|\rho_1|$.
Additionally, the increasing of the number of looks tends to correct this behaviour.

\section{Conclusion}

In this paper, 
we have derived 
four contrast measures 
in terms of
the correlation coefficient and of
the number of looks.
Using asymptotic results 
for theses measures, 
two methodologies based on the 
Kullback-Leibler and Hellinger distances  
were proposed 
as new discrimination techniques 
between PolSAR image regions.
These methods 
were applied to actual data 
and
the obtained results present evidence in favor
of both the proposals.

\bibliographystyle{IEEEtran}
\bibliography{../bibtexart}

\end{document}